\documentclass[a4paper, 11pt]{amsart}

\usepackage{amsfonts,amsmath,amssymb,amsxtra,amsthm,tikz}
\usepackage{mathrsfs}
\usepackage{times}
\usepackage{anysize}
\usepackage{graphicx}
\marginsize{1in}{1in}{1in}{1in}

\input xy
\xyoption{all}
\usepackage[all,knot]{xy}

\newcommand{\MQ}{\operatorname{\bf{Q}}}
\newcommand{\MUU}{\operatorname{\bf{U}}}

\newcommand{\Sub}{\operatorname{Sub}}
\newcommand{\Lk}{\operatorname{Lk}}
\newcommand{\St}{\operatorname{St}}
\newcommand{\mc}{\operatorname{\bf{c}}}
\newcommand{\mw}{\operatorname{\bf{w}}}
\newcommand{\ms}{\operatorname{\bf{s}}}
\newcommand{\mq}{\operatorname{\bf{q}}}
\newcommand{\ma}{\operatorname{\bf{a}}}
\newcommand{\me}{\operatorname{\bf{e}}}
\newcommand{\mo}{\operatorname{\bf{o}}}

\definecolor{cof}{RGB}{219,144,71}
\definecolor{pur}{RGB}{186,146,162}
\definecolor{greeo}{RGB}{91,173,69}
\definecolor{greet}{RGB}{52,111,72}
  \definecolor{lgreen} {RGB}{180,210,100}
  \definecolor{dblue}  {RGB}{119,158,201}
  \definecolor{lred}   {RGB}{220,0,0}
  \definecolor{nred}   {RGB}{224,0,0}
  \definecolor{norange}{RGB}{230,120,20}
  \definecolor{nyellow}{RGB}{255,221,0}
  \definecolor{ngreen} {RGB}{34,139,50}
  \definecolor{dgreen} {RGB}{78,138,21}
  \definecolor{nblue}  {RGB}{28,130,185}
 \definecolor{jblue}  {RGB}{0,103,166}

\theoremstyle{plain}

\newtheorem{Lm}{Lemma}
\newtheorem{Th}{Theorem}
\newtheorem{Cor}{Corollary}

\theoremstyle{definition}

\newtheorem{Ex}{Example}

\newenvironment{Pf}
{\noindent{\bf Proof\/}.}{{ $\Box$}\smallskip\par}

\title{Subword complexes and nil-Hecke moves}
\author{Mikhail Gorsky}

\address{\small Steklov Mathematical Institute, Gubkina str., 8, Moscow, 119991, Russia; \newline Universit\'e Paris Diderot -- Paris 7, Institut de Math\'ematiques de Jussieu -- Paris Rive Gauche,
B\^at. Sophie Germain, 75205 Paris Cedex 13, France}
\email{gorskym@math.jussieu.fr}

\begin{document}

\begin{abstract}
For a finite Coxeter group $W,$ a subword complex is a simplicial complex associated with a pair $(\MQ, \rho),$ where $\MQ$ is a word in the alphabet of simple reflections, $\rho$ is a group element. We describe the transformations of such a complex induced by nil-moves and inverse operations on $\MQ$ in the nil-Hecke monoid corresponding to $W.$ If the complex is polytopal, we also describe such transformations for the dual polytope. For $W$ simply-laced, these descriptions and results of \cite{Go} provide an algorithm for the construction of the subword complex corresponding to $(\MQ, \rho)$ from the one corresponding to $(\delta(\MQ), \rho),$ for any sequence of elementary moves reducing the word $\MQ$ to its Demazure product $\delta(\MQ).$ The former complex is spherical if and only if the latter one is empty.
\end{abstract}
\maketitle

\section {Introduction}


Subword complexes were introduced by A.~Knutson and E.~Miller in \cite{KM1} in the context of Schubert polynomials and matrix Schubert varieties. It was soon realized that they are interesting from the point of view of Coxeter combinatorics. A subword complex $\Delta(\MQ, \rho)$ is associated with a pair $(\MQ, \rho),$ where $\MQ$ is a word in the alphabet of simple reflections and $\rho$ is an element of a Coxeter group $W.$ The simplices in the complex correspond to the subwords in $\MQ$ whose complements contain reduced expressions of $\rho.$ The exchange axiom arises as the transition between two adjacent maximal simplices. In \cite{KM1}, subword complexes were shown to be vertex-decomposable and, therefore, shellable. This provides a new proof (and a new interpretation) of the Cohen-Macaulayness for matrix Schubert varieties and also for ordinary Schubert varieties, cf. \cite{KM1}. Using the shellability, in \cite{KM2}, Knutson and Miller proved that an arbitrary subword complex is homeomorphic to a sphere or to a ball. They also showed that the complex is spherical if and only if the Demazure product $\delta(\MQ)$ of $\MQ$ equals $\rho.$ This product can be defined as the unique reduced expression of $\MQ$ in the 0-Hecke monoid corresponding to the group $W.$ 

For a spherical subword complex it is natural to ask, whether it is polar dual to some simple polytope. The answer to the general question is still unknown, but in some cases such polytopes are constructed. One of the constructions is due to V.~Pilaud and C.~Stump \cite{PS} whose polytopes are called {\it brick}. As polytopes dual to the subword complexes there arise such classical polytopes as generalized associahedra and certain cyclic polytopes, cf. \cite{CLS}\cite{PS}. For a certain class of words, subword complexes turn out to be isomorphic to cluster complexes introduced by S.~Fomin and A.~Zelevinsky \cite{FZ}. Moreover, some important objects from the theory of cluster algebras, such as $c-$vectors, can be interpreted via natural geometric realizations of these subword complexes and their fans \cite{CLS}. One class of subword complexes is called multi-cluster complexes. A number of properties of the category of representations of the corresponding quiver have a natural interpretation in terms of combinatorics of these complexes, cf. \cite{CLS},\cite{CP}. 

In \cite{Go}, the author described the transformations of subword complexes and their dual polytopes induced by braid moves of a word in the corresponding Coxeter group. Braid moves in the corresponding 0-Hecke monoid are the same operations. Another important class of elementary moves in 0-Hecke monoids are the nil-moves and inverse operations. A main result of this article is the description of the transformations of subword complexes and their dual polytopes induced by nil-Hecke moves and inverse operations on a word. It follows from the classical Word Property of Coxeter groups, that any word in a 0-Hecke monoid can be reduced by a sequence of braid and nil-Hecke moves. Thus, every subword complex $\Delta(\MQ, \rho)$ can be obtained from a complex $\Delta(\delta(\MQ), \rho)$ by the transformations mentioned above, and any sequence of elementary moves connecting $\delta(\MQ)$ with $\MQ$ provides a sequence of such elementary transformations. In particular, if $\delta(\MQ) = \rho,$ the complex $\Delta(\delta(\MQ), \rho)$ is just an empty complex. We hope that this might be an important step on the way towards the answer to the general question of polytopality of subword complexes.

When $W$ is simply-laced, these elementary transformations are compositions of edge subdivions, inverse operations and taking the suspension. Edge subdivisions and $2-$truncations attracted the attention recently because of their applications to enumerative polynomials, in particular, to Gal's conjecture \cite{Ga} claiming the non-negativity of the coefficients of the $\gamma-$polynomial of an arbitrary flag spherical simplicial complex; details can be found in \cite{BV}. The class of polytopes which can be obtained from the cube of a fixed dimension by a sequence of $2-$truncations turns out to be very interesting. In particular, it contains all flag nestohedra, graph-associahedra and generalized associahedra of types $ABCD;$ details can be found in \cite{BV}. There is also an old and little-known theorem of M.H.A. Newman\cite{N} stating that any two of PL-homeomorphic complexes can be related to each other by a sequence of edge subdivisions and inverse operations. 
 
The paper is organized as follows. In Section 2 we recall the notion of edge subdivisions of simplicial complexes and of $2-$truncations of simple polytopes. Section 3 and Section 4 contain the introductions to finite Coxeter systems and to corresponding 0-Hecke monoids, respectively. Section 5 contains the definition and examples of subword complexes. In Section 6, we formulate and prove all new results.

This is a part of my ongoing Ph.D. project at Steklov Mathematical Institute. I am very grateful to my doctoral adviser Prof. Victor M. Buchstaber for the inspiration in the work and for his support and patience. Many thanks to Jean-Philippe Labb\'e, Vincent Pilaud and Salvatore Stella for explaining me the nature of the subword complexes and to Sergei Fomin who pointed out to me the relevance of the nil-Hecke monoids in this context. I am grateful to Alexander A. Gaifullin for informing me about the main result of the article \cite{N}. The work was supported by DIM RDM-IdF of the R\'{e}gion \^{I}le-de-France.

\section{Edge subdivisions and 2-truncations}
Given a simplex $\sigma$ in a simplicial complex $X,$ the {\it link} and the {\it star} of $\sigma$ are the following subcomplexes:
$$\Lk_X(\sigma) = \left\{\rho \in X | \sigma \cup \rho \in X, \sigma \cap \rho = \emptyset\right\};$$
$$\St_X(\sigma) = \left\{\rho \in X | \sigma \cup \rho \in X\right\}.$$

A {\it join} of two complexes $X_1$ and $X_2$ is a complex on the disjoint union of their vertex sets defined as follows:
$X_1 * X_2 = \left\{\sigma_1 \cup \sigma_2 | \sigma_1 \in X_1, \sigma_2 \in X_2\right\}.$
A {\it suspension} of a complex is its join with a disjoint union of two points.

Let $X$ be a simplicial complex. Let $\eta = \left\{s, t\right\}$ be an edge. Define $\Sub_{\eta}(X)$ to be
a simplicial complex constructed from $X$ by bisection of all simplices containing $\eta.$ In other words,
let $r$ be any letter not in the vertex set of $X.$ Then $S \cup \left\{r\right\}$ is the vertex set of $\Sub_{\eta}(X)$ and
$$\Sub_{\eta}(X) = \left\{\sigma|\eta \not\subset \sigma \in X\right\} \cup \left\{\sigma \cup \left\{r\right\}, \sigma \cup \left\{r, s\right\}, \sigma \cup \left\{r, t\right\} | \sigma \in \Lk_X(\eta)\right\} = $$
$$= \left\{\sigma|\eta \not\subset \sigma \in X\right\} \cup \left\{\sigma \cup \left\{r\right\} | \sigma \in \partial \St_X(\eta) \right\}.$$
$\Sub_{\eta}(X)$ is called a subdivision of $X$ along $\eta.$

A convex polytope of dimension $n$ is said to be {\it simple} if each of its vertices is contained in precisely $n$ facets. A polytope $P$ is simple if and only if its polar dual polytope $P^*$ is {\it simplicial}, i.e. each face of $P^*$ is a simplex.
The boundary of a simplicial $n-$dimensional polytope is a simplicial complex
of dimension $(n - 1).$ For a simple polytope $P,$ we shall denote by $K_P$ the
boundary complex $\partial P^*$ of the dual polytope. 
We say that $K_P$ is the {\it nerve complex} of $P.$ 
For simple polytopes, there is a dual operation to the edge subdivision: 
this is the truncation of a face of codimension $2,$ or simply the $2-${\it truncation}. 

Let $G$ be a face of codimension $2$ of a simple (combinatorial) polytope $P.$ Let $K_P = \partial P^*$ be its nerve complex and $J \in P^*$ be the face dual to $G.$ We say that the polytope $\widetilde{P},$ such that $K_{\widetilde{P}} = \Sub_J(K_P),$ is the {\it truncation} of $P$ at $G.$ Such a polytope $\widetilde{P}$ exists and is unique, up to combinatorial isomorphism.

Assume that $P$ is $n-$dimensional. Geometrically, the polytope $\widetilde{P}$ can be obtained from a realization of $P$ by intersecting the latter with a new half-space $H,$ such that the intersection of the $(n-1)$-dimensional plane $h = \partial H$ with $\partial P$ is precisely the link of $G$ in $\partial P.$ It means that the facets $\widetilde{P}$ are precisely the facets of $P$ and there is a one new facet $K$ isomorphic to $I \times G,$ and such that $\Lk_{\widetilde{P}}(K) = \Lk_P(G).$ Here $I$ is the closed interval $[0,1].$ More detailed treatment of $2-$truncations can be found in \cite{BV}.

\section{Coxeter groups}
Let $W$ be a finite Coxeter group, that is, a finite group generated by a set $S = \left\{s_1,\ldots,s_n\right\}$ of {\it simple reflections}, modulo the relations $(s_i s_j)^{m_{ij}} = e.$ Here 
$m_{ii} = 2, m_{ij} = m_{ji}.$ 

$W$ is \textit{{simply-laced}} if $m_{ij} \in \left\{2,3\right\},$ for any $i,j.$ The \textit{{length}} $l(w)$ of an element $w \in W$ is the length of the smallest
expression of $w$ as a product of the generators in $S.$ An expression $w = w_1 w_2 \ldots w_p$
with $w_1,\ldots,w_p \in S$ is called \textit{{reduced}} if $p = l(w).$ We denote by $w_o$ the longest element in $W,$ it is known to be unique.
If the word $\mw$ in the alphabet $S$ contains a subword $\ms_i\ms_j\ms_i\ms_j\ms_i\ldots$ of length $m_{ij},$ then there is a {\it braid move} transforming $\mw$ into $\mw'$ by changing $\ms_i\ms_j\ms_i\ms_j\ms_i\ldots$ by the subword $\ms_j\ms_i\ms_j\ms_i\ms_j\ldots$ of the same length $m_{ij}.$

We denote by $S^*$ the set of words on the alphabet $S,$ and by $\me$ the empty word. 
To avoid confusion, we denote with a square letter $\ms$ the letter of the alphabet $S$
corresponding to the single reflection $s \in S.$ Similarly, we use a square letter like $\mw$
to denote a word of $S^*,$ and a normal letter like $w$ to denote its corresponding group
element in $W.$ 

On $S^*,$ there are two types of operations reflecting the group structure of $W.$ A {\it nil-Coxeter move} removes two consecutive identical letters $\ms_i\ms_i$ from a word $\mw \in S^*,$ for some $i.$ If $\mw$ contains a subword $\ms_i\ms_j\ms_i\ms_j\ms_i\ldots$ of length $m_{ij},$ then there is a {\it braid-move} transforming $\mw$ into $\mw'$ by changing $\ms_i\ms_j\ms_i\ms_j\ms_i\ldots$ by the subword $\ms_j\ms_i\ms_j\ms_i\ms_j\ldots$ of the same length $m_{ij}.$ Note that neither a nil-move, nor a braid-move changes a group element $w \in W$ expressed by $\mw.$ Recall the {\bf Word Property} which holds for any Coxeter system $(W; S).$

\begin{Th} [{\cite[Theorem~3.3.1]{BB}}] \label{wordproperty}
Any expression $\mw$ for $w \in W$ can be transformed into a
reduced expression for $w$ by a sequence of Coxeter nil-moves and braid-
moves.
\end{Th}

A {\it Coxeter element} $c$ is a product of
all simple reflections in some order. We choose an arbitrary
reduced expression $\mc$ of $c$ and denote by $\mw(\mc)$ the {\it $\mc$-sorting word} of $w,$ that is the
lexicographically first (as a sequence of positions) reduced subword of ${\mc}^{\infty} = \mc\mc\mc\ldots$ for $w.$ In
particular, ${\mw}_{\mo}(\mc)$ denotes the $\mc$-sorting word of the longest element $w_o \in W.$

\section{0-Hecke monoid}

For a finite Coxeter group $W$ with
generators $S,$ there is defined a corresponding {\it $0-$Hecke monoid}. 
It has a generating set $X = {x_1,\ldots, x_n},$ where $x_i$ corresponds to $s_i.$ The difference is that each Hecke generator is an idempotent, namely satisfies
$x_i^2 = x_i,$ whereas the Coxeter generators are involutions. In other words, instead of reflections one has projectors.  Also,
each Coxeter braid relation $(s_is_j)^{m_{ij}} = e$ provides a corresponding
Hecke braid relation of the form 
$$x_ix_jx_i\ldots = x_jx_ix_j \ldots ,$$
with $m_{ij}$ alternating terms on each side. 

The Word Property (Theorem \ref{wordproperty}) implies that an arbitrary word in the $0-$Hecke generators may be
reduced by a sequence of Hecke nil-moves $x^2_i \to x_i$ and braid moves $x_ix_jx_i\ldots \to 
x_jx_ix_j \ldots.$ The element of the monoid expressed by the result does not depend on the sequence of moves; this is precisely the product of letters of the initial word in the monoid. A reduced word in the $0-$Hecke generators corresponds to a
reduced word in the Coxeter generators by switching occurrences of $x_i$
to $s_i.$ For any word $\MQ$ in the alphabet $S,$ we define the {\it Demazure product} (or {\it $0-$Hecke product}) $\delta(\MQ)$ as follows: we switch all letters $s_i$ to $x_i,$ reduce the result, switch all $x_i$ back to $s_i,$ and consider the result as an element of $W.$ The above arguments show that this definition is correct; in other words, the $0-$Hecke product is a multiplication, inducing a monoid structure on $W.$  It is known and easy to check that $\MQ$ contains a reduced expression for some group element $\rho \in W$ if and only if $\delta(\MQ)$ does. The $0-$Hecke product is actually the multiplication rule for basis elements in the $0-$Hecke algebra attached to $W.$ An overview of its properties may be found in \cite{BM}.

\section{Subword complexes}

Let $W$ be a finite Coxeter group with the set of simple reflections $S,$ let $\MQ := \mq_1\ldots\mq_m$ belong to $S^*$
and let $\rho$ be an element in $W.$ The {\it subword
complex} $\Delta(\MQ; \rho)$ is the pure simplicial complex of subwords of $\MQ,$ whose complements
contain a reduced expression of $\rho.$ The vertices of this simplicial complex are labeled
by (positions of) the letters in the word $\MQ.$ Note that two positions are different
even if the letters of $\MQ$ at these positions coincide. 
The maximal simplices of the subword complex $\Delta(\MQ; \rho)$ are the complements of reduced expressions of $\rho$ in the word $\MQ.$

In \cite{KM1}, it was shown that the subword complex $\Delta(\MQ; \rho)$ is either a triangulated sphere (or simply {\it spherical}), or a triangulated ball. It is spherical if and only if the Demazure product $\delta(\MQ)$ is equal to $\rho,$ see the proof in \cite[Section~3]{KM1}.  In some generality, these spherical subword complexes are polar dual to certain simple polytopes. The general description of these polytopes and their geometric realizations is not yet known; under certain conditions, such realizations are constructed and called the {\it brick polytopes} in \cite{PS}. We are interested in combinatorial types of these polytopes. The facets of such a polytope are labeled by (positions of) the letters in the word $\MQ,$ whose complements contain a reduced expression of $\rho.$ Let us give some examples. Denote for the Coxeter group $A_n$ the simple reflections by $s_1,\ldots,s_n$ in a natural order. 

\begin{Ex} \label{cluster} (Cluster complexes and generalized associahedra).
For any $\mbox{c},$ the subword complex \newline $\Delta(\mc{\mw}_o(\mc); w_o)$ coincides with the cluster complex of type $W.$ Its dual polytope is a generalized associahedron of type $W.$ For example, for $W$ of type $A_n$ and
$\mc = \ms_1\ms_2\ldots\ms_n,$ we have
$$\mc{\mw}_o(\mc) = \ms_1\ms_2\ldots\ms_n\ms_1\ms_2\ldots\ms_n\ms_1\ms_2\ldots{\ms}_{n-1}\ms_1\ms_2\ldots\ms_{n-2}\ldots\ms_1\ms_2\ms_1.$$
Every position yields a vertex of the complex and the facet of the polytope. 
\end{Ex}

\begin{Ex} (Multi-cluster complexes).
For any $\mc,$ the subword complex $\mathcal{S}({\mc}^k\mw_o(\mc))$ is called a {\it $k$-cluster complex} of type $W.$ It is not known in general, whether thus defined complexes are polar dual to some polytopes. When this holds, the polytope polar dual to $\mathcal{S}(\mbox{c}^k\mbox{w}_o(c))$ is called the {\it $k-$associahedron} of type $W.$ See {\cite[PS]{CLS}} for more details.
\end{Ex}

We will use the following simple observation.

\begin{Lm} [\cite{KM2}] \label{linksubword}
The link of a simplex corresponding to a word $\ma_1 \ma_2\ldots\ma_x$ in a subword complex $\Delta(\MUU; \rho),$ where $\MUU \in S^*, \rho \in W$ is isomorphic to the complex $\Delta(\MUU'; \rho),$ where $\MUU'$ is obtained from $\MUU$ by removing all the letters $\ma_i, i = 1,2,\ldots,x.$
\end{Lm}

\section{Geometry of nil-moves}

In this section, we describe the relation between subword complexes for words, linked by a nil-move.

\begin{Th} \label{nil}
Assume that ${\MQ}^{'}$ is a word obtained from $\MQ$ by doubling one letter: 
$$\ldots \mq \ldots \quad \longrightarrow \quad \ldots \mq\mq \ldots$$
Then, for any $\rho \in W,$ the following holds:
\begin{itemize}
\item $\Delta({\MQ}^{'}, \rho)$ either coincides with the suspension $\Sigma(\Delta(\MQ, \rho)),$ or can be obtained from it by an inverse edge subdivision. $\Delta(\MQ, \rho)$ is the link in $\Delta({\MQ}^{'}, \rho)$ of the vertex corresponding to one of the doubling letters $\mq.$

\item $\Delta({\MQ}^{'}, \rho)$ is polytopal iff $\Delta(\MQ, \rho)$ is. In such a case, $B({\MQ}^{'}, \rho)$ either coincides with $B(\MQ, \rho) \times I,$ or can be obtained from it by an inverse $2-$truncation. $B(\MQ, \rho)$ is the facet of $B({\MQ}^{'}, \rho)$ corresponding to one of the doubling letters $\mq.$
\end{itemize}
\end{Th}

\begin{Pf}
Denote doubling positions of the letter $\mq$ in ${\MQ}^{'}$ by ${\mq}^1$ and ${\mq}^2,$ its position in $\MQ$ by ${\mq}^0.$ 
No reduced expression of $\rho$ can contain two consecutive identical letters; thus, the complement to every reduced expression of $\rho$ in ${\MQ}^{'}$ contains at least one of ${\mq}^1$ and ${\mq}^2.$ Therefore, any maximal simplex in $\Delta({\MQ}^{'}, \rho)$ is contained in the union $(\Lk_{\Delta({\MQ}^{'};\rho)}(\left\{{\mq}^{1}\right\}) * \left\{{\mq}^{2}\right\}) \cup (\Lk_{\Delta({\MQ}^{'};\rho)}(\left\{{\mq}^{2}\right\}) * \left\{{\mq}^{1}\right\}).$
Since each simplex in $\Delta({\MQ}^{'}; \rho)$ is contained in a maximal one, we have 
$$\Delta({\MQ}^{'}; \rho) \subset (\Lk_{\Delta({\MQ}^{'};\rho)}(\left\{{\mq}^{1}\right\}) * \left\{{\mq}^{2}\right\}) \cup (\Lk_{\Delta({\MQ}^{'};\rho)}(\left\{{\mq}^{2}\right\}) * \left\{{\mq}^{1}\right\}).$$
The inverse inclusion is clear; thus, we have
$$\Delta({\MQ}^{'}; \rho) = (\Lk_{\Delta({\MQ}^{'};\rho)}(\left\{{\mq}^{1}\right\}) * \left\{{\mq}^{2}\right\}) \cup (\Lk_{\Delta({\MQ}^{'};\rho)}(\left\{{\mq}^{2}\right\}) * \left\{{\mq}^{1}\right\}).$$
By Lemma \ref{linksubword}, we have 
$$\Lk_{\Delta({\MQ}^{'}, \rho)}({\mq}^1)  \cong \Delta(\MQ, \rho) \cong \Lk_{\Delta({\MQ}^{'}, \rho)}({\mq}^2),$$
where the first isomorphism is given by switching from ${\mq}^2$ to ${\mq}^0,$ the second one is given by switching from ${\mq}^0$ to ${\mq}^1.$ 
We see that $\Delta({\MQ}^{'}, \rho)$ equals to the following union:
$$\Delta({\MQ}^{'}, \rho) = \left\{\sigma, \sigma \cup \left\{{\mq}^{1}\right\}, \sigma \cup \left\{{\mq}^{2}\right\} | \sigma \in (\Delta(\MQ, \rho) \backslash \St_{\Delta(\MQ, \rho)}({\mq}^0)) \right\} \sqcup $$
$$ \sqcup (\Lk_{\Delta(\MQ, \rho)}({\mq}^0)) *  (\left\{{\mq}^1, {\mq}^2\right\}).$$

First assume that ${\mq}^0$ is not a vertex of the complex $\Delta(\MQ, \rho).$ Then the second set in the big union is empty, and we have 
$$\Sigma(\Delta({\MQ}^{'}, \rho)) = \Sigma(\Delta(\MQ, \rho) \backslash \St_{(\Delta(\MQ, \rho))}({\mq}^0)) = \Sigma(\Delta(\MQ, \rho)).$$
Assume now that ${\mq}^0$ is a vertex of the complex $\Delta(\MQ, \rho).$ Then it is easy to check that  the suspension $\Sigma(\St_{\Delta(\MQ, \rho)}({\mq}^0))$ is the subdivision of the second term at the edge $\left\{{\mq}^1, {\mq}^2\right\},$ where the new vertex is ${\mq}^0.$ Thus, the second term is obtained from $\Sigma(\St_{\Delta(\MQ, \rho)}({\mq}^0))$ by an inverse edge subdivision, and the whole union $\Delta({\MQ}^{'}, \rho)$ is obtained from 
$$\Sigma(\Delta(\MQ, \rho) \backslash \St_{(\Delta(\MQ, \rho))}({\mq}^0)) \sqcup \Sigma(\St_{\Delta(\MQ, \rho)}({\mq}^0)) = \Sigma(\Delta(\MQ, \rho))$$
by an inverse edge subdivision (one eliminates the vertex ${\mq}^0$ to obtain an edge $\left\{{\mq}^1, {\mq}^2\right\}$).
This finishes the proof of the first statement. The proof of the second statement is dual, modulo the polytopality part, which is trivial.
\end{Pf}

\begin{Ex}

Take $W = A_2, \rho = w_o$ and $\MQ = \ms_1\ms_2\ms_1\ms_2\ms_1.$ As we know, $\Delta(\MQ, \rho)$ is a pentagon, and the dual polytope $P = As^2$ is again a pentagon. Then the doubling of the first letter in $\MQ$ transforms $P$ into $P'$ which is $P \times I$ with one truncated edge.

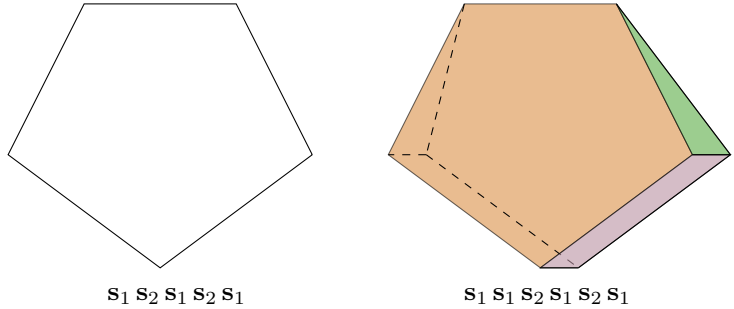
\begin{center}
\begin{tikzpicture}[scale=0.5]
\coordinate (B1) at (-1,7) {};
\coordinate (B3)  at (-9,7) {};
\coordinate (B2) at (-5,4) {};
\coordinate (B4)  at (-7,11) {};
\coordinate (B5)  at (-3,11) {};
\draw (B1) --  (B2) --  (B3) --  (B4) -- (B5) -- (B1);

\coordinate (A1) at (9,7) {};
\coordinate (A3)  at (1,7) {};
\coordinate (A2) at (5,4) {} {};
\coordinate (A4) at (3,11) {} {};
\coordinate (A5) at (7,11) {} {};
\draw [fill=cof,opacity=0.6]
(A1) -- (A2) --  (A3) --  (A4) --  (A5) --  (A1);
\coordinate (v1) at (2,7) {};
\coordinate (v2) at (6,4) {};
\coordinate (v3) at (10,7) {};
\draw[dashed] (A4) -- (v1);
\draw[dashed] (v1) -- (v2);
\draw  [fill=pur,opacity=0.6]
(v2) -- (v3) -- (A1) -- (A2) -- (v2);
\draw [fill=greeo,opacity=0.6]
(v3) -- (A1) -- (A5) -- (v3);
\draw (A5) edge (v3);
\draw (v3) edge (A1);
\draw[dashed] (v1) edge (A3);
\draw (v2) edge (A2);
\draw (v2) -- (v3);

\end{tikzpicture}

\small $\ms_1\ms_2\ms_1\ms_2\ms_1 \qquad\qquad\qquad\qquad \ms_1\ms_1\ms_2\ms_1\ms_2\ms_1$
\end{center}

Note that we can double any letter in $\MQ$ and get the same result. This holds since any letter in $\MQ$ gives a facet of $P,$ and all facets in the pentagon are combinatorially indistinguishable.
\end{Ex}

Recall the main result of \cite{Go}:

\begin{Th} [\cite{Go}] \label{braid}
Assume that words $\MQ$ and ${\MQ}^{'}$ are related by a braid move. Then, under certain conditions, $\Delta({\MQ}^{'}; \rho)$ can be obtained from $\Delta(\MQ; \rho)$ by a sequence of several subdivisions along one edge and inverse edge subdivisions along another one. This holds for any $\MQ, {\MQ}^{'},$ when $W$ is simply-laced. The dual result also holds for polytopes and $2-$truncations.
\end{Th}

\begin{Cor}
Suppose $W$ is simply-laced. Then, for any $\MQ, \rho,$ any sequence of braid-moves and inverse nil-Hecke moves from $\delta(\MQ)$ to $\MQ$ provides a sequence of elementary transformations  described in Theorems \ref{nil} and \ref{braid} from $\Delta(\delta(\MQ); \rho)$ to $\Delta(\MQ; \rho).$ At least one such sequence exists (by the Word Property). If $\delta(\MQ) = \rho,$ there exists a sequence of elementary transformations from the empty complex to $\Delta(\MQ; \rho).$ If, moreover, $\Delta(\MQ; \rho)$ is polytopal, its dual polytope can be obtained from the empty complex by a sequence of elementary transformations.
\end{Cor}

\renewcommand{\refname}{References}

\medskip

\end{document}